\documentclass[11pt,a4paper]{amsart}
\usepackage[T1]{fontenc}
\usepackage{lmodern}
\usepackage{amsmath,amssymb,amsthm,mathtools}
\usepackage{geometry,microtype,xcolor,hyperref}
\hypersetup{colorlinks=true,linkcolor=blue!45!black,citecolor=blue!45!black,urlcolor=blue!50!black,pdftitle={Non-migration for planar elastic flows},pdfauthor={Matteo Novaga}}
\numberwithin{equation}{section}
\newtheorem{theorem}{Theorem}[section]
\newtheorem{proposition}[theorem]{Proposition}
\newtheorem{lemma}[theorem]{Lemma}
\newtheorem{corollary}[theorem]{Corollary}
\theoremstyle{remark}
\newtheorem{remark}[theorem]{Remark}
\newcommand{\R}{\mathbb R}
\newcommand{\T}{\mathbb T}
\newcommand{\norm}[1]{\lVert #1\rVert}
\author{Matteo Novaga}
\address{Dipartimento di Matematica, Universit\`a di Pisa, Pisa, Italy}
\email{matteo.novaga@unipi.it}
\title{Non-migration for planar elastic flows}
\date{\today}

\begin{document}
\begin{abstract}
We consider a closed planar curve evolving by the free or the length-preserving
elastic flow. We prove that the scale-invariant bound
$L\int k^2\,ds\le 4\pi^2+8/3$ at the initial time prevents the entire curve
from crossing any line between two times. For the free flow,
the proof combines a monotonicity estimate for the normalized bending energy
with a barrier that compares centroid motion to an expanding energy radius.
For the length-preserving flow, a quantitative gradient estimate controls the
total centroid displacement.
\end{abstract}

\maketitle

\section{Introduction}
Let $\gamma:\T\to\R^2$, $\T=\R/\mathbb Z$, be a smooth closed immersed
curve. We denote by $s$ the arclength parameter, by $T=\partial_s\gamma$ the unit
tangent, and by $\nu$ the left unit normal (the counterclockwise rotation of
$T$ by $\pi/2$), so that $T_s=k\nu$. Euler's bending energy, the length,
and the normalized bending energy of $\gamma$ are
\begin{equation}\label{eq:def-U}
 B(\gamma):=\int_\gamma k^2\,ds,\qquad
 L(\gamma):=\int_\gamma ds,\qquad U(\gamma):=L(\gamma)B(\gamma).
\end{equation}
Note that the product $U$ is invariant under dilations.

We consider the $L^2(ds)$-gradient flow of $B$ and its projection onto the
fixed-length constraint:
\begin{equation}\label{eq:flows}
 \partial_t\gamma=-G\nu,\qquad
 G=\begin{cases}
 G_0:=2k_{ss}+k^3, & \text{free flow},\\
 G_c:=2k_{ss}+k^3-\lambda_c(t)k,
 & \text{length-preserving flow},
 \end{cases}
\end{equation}
where
\begin{equation}\label{eq:lambda-intro}
 \lambda_c(t)=\frac{\int_{\gamma_t}(k^4-2k_s^2)\,ds}
 {\int_{\gamma_t}k^2\,ds}.
\end{equation}
Thus $\lambda_c$ is chosen so that $\int_{\gamma_t}kG_c\,ds=0$. We write the
evolution in normal gauge; tangential reparametrizations do not affect the
geometric statements below. 
By~\cite{DKS}, both flows exist for all times
starting from any smooth closed curve. 

The study of elastic flows of closed curves goes back at least to
Wen~\cite{Wen} and Polden~\cite{Polden}. Dziuk--Kuwert--Sch\"atzle~\cite{DKS}
established global existence for the free, length-penalized, and
length-preserving flows of smooth closed curves. The survey of
Mantegazza--Pluda--Pozzetta~\cite{MPPSurvey} gives an account of the
analysis of elastic flows of curves and networks; see also
Barrett--Garcke--N\"urnberg~\cite{BarrettGarckeNuernberg} for numerical
approximations. For open curves, Dall'Acqua--Lin--Pozzi studied the
fixed-length flow with natural boundary conditions~\cite{DallAcquaLinPozziNatural};
clamped-end versions were treated in~\cite{DallAcquaLinPozzi,RuppSpener},
and a free-boundary version in~\cite{DallAcquaSchlierf}.
In the free closed case, circles expand rather than approach a stationary
curve. Near-circle stability after rescaling was proved first under a
curvature-derivative condition by Miura--Wheeler~\cite{MiuraWheeler} and
then under a small curvature-oscillation condition by
Andrews--Wheeler~\cite{AndrewsWheeler}.

A question raised by Huisken asks whether a closed planar curve under the
length-penalized elastic flow can move entirely from one side of a line to
the other~\cite{MPPSurvey,KemmochiMiura}. Such \emph{migration} has been
constructed for open curves with natural boundary conditions, both for the
length-preserving flow~\cite{KemmochiMiura} and for the length-penalized elastic
flow with constant multiplier~\cite{MiuraII}. These examples do not resolve
the question for closed curves. Miura--Wheeler~\cite{MiuraWheeler} also identify the migration
question for the free closed flow as open. Our result supplies a quantitative
non-migration criterion for the two closed-curve flows in~\eqref{eq:flows}.

Fenchel's theorem and Cauchy--Schwarz give $U\ge4\pi^2$ for every closed
immersed curve, with equality exactly for a once-covered circle. If $n$ is
the turning number of $\gamma$, then
\begin{equation}\label{eq:oscillation}
 L\int_\gamma\left(k-\frac{2\pi n}{L}\right)^2\,ds
 =U-4\pi^2n^2.
\end{equation}
The defect on the right is the normalized curvature oscillation used
in~\cite{AndrewsWheeler}. Our assumption concerns $U$ alone: it does not
impose a turning number, embeddedness, or smallness of curvature derivatives.

\begin{theorem}[Non-migration result]\label{th:main}
Let $\gamma:\T\times[0,+\infty)\to\R^2$ be a smooth family of closed regular immersed curves solving either equation in~\eqref{eq:flows}. Suppose that
\begin{equation}\label{eq:main-hyp}
 U(\gamma_0)=L(\gamma_0)\int_{\gamma_0}k^2\,ds
 \le4\pi^2+\frac83.
\end{equation}
Then, for all $0\le t_0<t_1$, we have
\begin{equation}\label{eq:main-persist}
 4\pi^2\le U(\gamma_{t_1})\le U(\gamma_{t_0})\le U(\gamma_0).
\end{equation}
In particular, every $\gamma_t$ is embedded. Moreover, for every unit
vector $e\in S^1$ and $h\in\R$, the following two inclusions cannot
both hold:
\begin{equation}\label{eq:halfplanes}
 \gamma_{t_0}\subset\{x:x\cdot e\ge h\},\qquad
 \gamma_{t_1}\subset\{x:x\cdot e\le h\}.
\end{equation}
\end{theorem}

The energy threshold in Theorem~\ref{th:main} is sufficient, and we do not
claim that it is optimal. Note that $4\pi^2+8/3\simeq42.15$ lies well below
Miura's sharp threshold $4\varpi^*\simeq112.44$: every closed immersed curve
with $U<4\varpi^*$ is embedded, and the figure-eight elastica attains
equality~\cite{MiuraLiYau}. For the length-preserving flow, this threshold
also yields embeddedness at all times~\cite[Theorem~4.2]{MiuraLiYau}, but
embeddedness does not by itself prevent the curve from crossing a line. In
our smaller sublevel, embeddedness follows from an elementary argument
(Lemma~\ref{lem:topology}). The non-migration conclusion is stronger and
requires control of the position of the curve.

The geometric argument uses the length centroid
\begin{equation}\label{eq:cent-R}
 c=\frac1L\int_\gamma\gamma\,ds
\end{equation}
and the directional margins
\begin{equation}\label{eq:margins}
 \rho_e^+=\max_\gamma\gamma\cdot e-c\cdot e,
 \qquad\rho_e^-=c\cdot e-\min_\gamma\gamma\cdot e.
\end{equation}
For the free flow, $B$ decreases strictly and the auxiliary energy radius
$R_E=2\pi/B$ expands. Dissipation of $B$ alone does not imply monotonicity
of $U=LB$, since $L$ varies. We first prove that the sublevel set defined
by~\eqref{eq:main-hyp} is preserved and that $U$ is nonincreasing on it; then
we establish
\begin{equation}\label{eq:barrier-intro}
 |c'|\le\frac34R_E',\qquad
 \rho_e^\pm\ge\frac34R_E\quad(e\in S^1).
\end{equation}
For the length-preserving flow, $L$ is constant and a quantitative gradient
estimate shows that the total displacement of $c$ between any two times is
less than $L/(8\pi^2)$; the directional margins are at least $3L/(8\pi)$.
Both estimates exclude a complete crossing. We establish the variation
formulas and the topological facts in Section~\ref{sec:prelim}, the
free-flow energy estimate in Section~\ref{sec:free}, the geometric bounds in
Section~\ref{sec:geom}, and the length-preserving estimate in
Section~\ref{sec:lp}; they are combined in Section~\ref{sec:proof}.

\section{Preliminaries}\label{sec:prelim}

\subsection{Variation formulas}
All integrals are with respect to arclength unless otherwise indicated.  Put
\begin{equation}\label{eq:Y-S}
 Y=\int_\gamma G^2\,ds,\qquad S=\int_\gamma(k^4-2k_s^2)\,ds.
\end{equation}
\begin{lemma}\label{lem:exact}
Along the free flow in~\eqref{eq:flows},
\begin{equation}\label{eq:basic}
 B'=-Y,\qquad L'=S,\qquad U'=BS-LY,
 \qquad R_E'=\frac{2\pi Y}{B^2}>0.
\end{equation}
Furthermore,
\begin{equation}\label{eq:translation}
 Lc'=\int_\gamma G\,[k(\gamma-c)+\nu]ds,
\end{equation}
\begin{equation}\label{eq:dilation}
 B=-\int_\gamma G\,(\gamma-c)\cdot\nu\,ds,
 \qquad\sqrt Y\ge\frac{2\pi B}{L^{3/2}}.
\end{equation}
\end{lemma}
\begin{proof}
For an arbitrary normal variation $\partial_t\gamma=f\nu$, the variation formulas corresponding to $T_s=k\nu$ are
\[
 \partial_tds=-kf\,ds,\qquad k_t=f_{ss}+k^2f.
\]
Taking $f=-G_0$, differentiation and periodic integration by parts give
\[
 B'=\int_\gamma[-2kG_{ss}-k^3G]ds=-\int_\gamma G^2ds,
 \qquad L'=\int_\gamma kG\,ds=S.
\]
This yields the first three formulas in~\eqref{eq:basic}; the fourth follows from $R_E=2\pi/B$.  There is no stationary smooth closed free elastica: by spatial scaling $\gamma\mapsto(1+\varepsilon)\gamma$, the first variation of $B$ equals $-B\ne0$, so $G$ cannot vanish identically.  Thus $Y>0$.

Translational invariance implies $\int_\gamma G\nu\,ds=0$.  Differentiating $L^{-1}\int\gamma\,ds$ then gives
\[
 Lc'=\int_\gamma[-G\nu+kG(\gamma-c)]ds
      =\int_\gamma G[k(\gamma-c)+\nu]ds,
\]
where the last equality again uses $\int G\nu\,ds=0$.
The scaling variation gives the first identity in~\eqref{eq:dilation}.  The elementary estimate
\begin{equation}\label{eq:cent-bound}
 \norm{\gamma-c}_{L^\infty}\le\frac L4
\end{equation}
results from writing $\gamma(s)-c=L^{-1}\int_0^L(\gamma(s)-\gamma(r))\,dr$ and bounding each difference by the length of the shorter connecting arc.  The Wirtinger inequality for $L$-periodic functions with zero mean gives the $L^2$ bound
\[
 \norm{\gamma-c}_2\le\frac{L}{2\pi}\norm T_2
 =\frac{L^{3/2}}{2\pi}.
\]
Consequently, by Cauchy--Schwarz, $B\le\sqrt Y\norm{\gamma-c}_2\le L^{3/2}\sqrt Y/(2\pi)$, proving~\eqref{eq:dilation}.
\end{proof}

\subsection{Low-energy topology}
The turning number and embeddedness needed below follow directly from the
energy bound.
\begin{lemma}\label{lem:topology}
Let $\gamma$ be a smooth closed unit-speed immersed planar curve. If
\begin{equation}\label{eq:topo}
 LB-4\pi^2<\frac{8\pi^2}{3},
\end{equation}
then its turning number is either $+1$ or $-1$. If, in addition,
$LB-4\pi^2\le8/3$, then $\gamma$ is embedded.
\end{lemma}
\begin{proof}
Write $\theta=2\pi s/L$ and identify $\tau(\theta)=T(L\theta/(2\pi))$ with a complex number of modulus one.  Closure implies
\(\int_0^{2\pi}\tau\,d\theta=0\).  Expand
\[
 \tau(\theta)=\sum_{j\in\mathbb Z\setminus\{0\}}a_je^{ij\theta}.
\]
Parseval gives
\[
 \sum_j|a_j|^2=1,\qquad
 \sum_jj^2|a_j|^2=\frac{LB}{4\pi^2}.
\]
With $d=LB-4\pi^2$ we have
\[
 \sum_{|j|\ge2}(j^2-1)|a_j|^2=\frac{d}{4\pi^2}.
\]
The elementary telescoping sum $\sum_{|j|\ge2}(j^2-1)^{-1}=3/2$ therefore yields
\begin{equation}\label{eq:tangent-high}
 \left\|\tau-(a_1e^{i\theta}+a_{-1}e^{-i\theta})\right\|_\infty
 \le\sqrt{\frac{3d}{8\pi^2}}<1.
\end{equation}
As $|\tau|=1$, the first-harmonic curve $\tau_1(\theta)=a_1e^{i\theta}+a_{-1}e^{-i\theta}$ never passes through the origin; straight-line homotopy between $\tau$ and $\tau_1$ also avoids the origin.  If $|a_1|=|a_{-1}|$, the function $\tau_1$ vanishes somewhere, a contradiction.  If $|a_1|>|a_{-1}|$, it is homotopic in $\mathbb C\setminus\{0\}$ to $a_1e^{i\theta}$; otherwise it is homotopic to $a_{-1}e^{-i\theta}$.  Thus $\deg\tau=\deg\tau_1=\pm1$.

For the second assertion reverse the orientation if necessary. Write
\[
 T(s)=e^{i\alpha}e^{2\pi is/L}e^{i\varphi(s)},\qquad
 \varphi_s=k-\frac{2\pi}{L},\qquad \int_0^L\varphi\,ds=0.
\]
Fourier series and
Cauchy--Schwarz yield
\[
 \norm\varphi_\infty\le
 \left(\sum_{j\ne0}\frac{L^2}{4\pi^2j^2}\right)^{1/2}
 \left(\frac1L\int_0^L\varphi_s^2\,ds\right)^{1/2}
 =\sqrt{\frac{LB-4\pi^2}{12}}\le\frac{\sqrt2}{3}.
\]
Let $s\ne s'$ in $[0,L)$. Exchanging them if necessary, we may write
$s'=s+\ell$ modulo $L$ with $0<\ell\le L/2$. Comparing the displacement along
this arc with the corresponding chord of a circle of radius $L/(2\pi)$, and
using $|e^{i\varphi}-1|\le|\varphi|$ and $\sin x\ge2x/\pi$ for
$x\in[0,\pi/2]$, we get
\[
 |\gamma(s+\ell)-\gamma(s)|
 \ge\frac L\pi\sin\frac{\pi\ell}{L}
       -\ell\norm\varphi_\infty
 \ge\ell\left(\frac2\pi-\frac{\sqrt2}{3}\right)>0.
\]
Hence $\gamma(s)\ne\gamma(s')$, and the curve has no self-intersections.
\end{proof}

\section{Signed evolution for the free flow}\label{sec:free}
In this section we establish the analytic input to Theorem~\ref{th:main} independently of the later geometry.  The point is to retain the favorable signs of fourth- and sixth-order curvature terms.

\begin{proposition}\label{prop:persist}
For any smooth closed free elastic flow with turning number $\pm1$, if $0<U-4\pi^2\le8/3$, then $U'<0$.  More precisely, after instantaneous ambient rescaling to $L=2\pi$ and an orientation choice so that $\int k\,ds=2\pi$, writing $k=1+u$, $\int u\,ds=0$, one has
\begin{equation}\label{eq:diss-bound}
\frac{U'}{L}\le-\frac1{250}\norm{u_{ss}}_2^2<0.
\end{equation}
Consequently $U(t)$ is nonincreasing as long as $U-4\pi^2\le8/3$, and the sublevel is forward invariant.
\end{proposition}
\begin{proof}
\emph{Instantaneous normalization.}  Both $U$ and the sign of $U'$ are invariant under ambient dilations (with the flow time rescaled by the fourth power of the spatial factor).  Hence we may freeze the rescaling at one time and work with $L=2\pi$, $k=1+u$, $\int u=0$.  Set
\begin{equation}\label{eq:VQX}
 \begin{gathered}
 V=\int u^2,\quad Q=\int u_s^2,\quad X^2=\int u_{ss}^2,\\
 I_j=\int u^j\quad(3\le j\le6),\qquad
 J_1=\int uu_s^2,\quad J_2=\int u^2u_s^2.
 \end{gathered}
\end{equation}
Thus $d:=U-4\pi^2=2\pi V$, and $V\le4/(3\pi)<1/2$.

\emph{Exact identity.}  Since $G=2u_{ss}+(1+u)^3$, direct expansion and periodic integration by parts yield
\begin{align}
 B&=2\pi+V,\nonumber\\
 S&=2\pi+6V+4I_3+I_4-2Q,\label{eq:BS-exp}\\
 Y&=4X^2-12Q-24J_1-12J_2
       +2\pi+15V+20I_3+15I_4+6I_5+I_6.\nonumber
\end{align}
Substituting into $U'=BS-LY$ gives
\begin{equation}\label{eq:U-exact}
 \begin{split}
 \frac{U'}L={}&-4X^2+10Q-8V-16I_3-14I_4-6I_5-I_6\\
 &\quad+24J_1+12J_2
 +\frac V{2\pi}(6V+4I_3+I_4-2Q).
 \end{split}
\end{equation}
By $J_1=-\frac12\int u^2u_{ss}$ and $J_2=-\frac13\int u^3u_{ss}$, the following identity is \emph{exact}:
\begin{equation}\label{eq:signed-square}
 \begin{split}
 &-4X^2+24J_1+12J_2-14I_4-6I_5-I_6\\
 &\qquad=-4\int\left(u_{ss}+\frac32u^2+\frac12u^3\right)^2ds-5I_4.
 \end{split}
\end{equation}
It follows that
\begin{equation}\label{eq:signed-evolution}
\begin{aligned}
 \frac{U'}L={}&-4\left\|u_{ss}+\frac32u^2+\frac12u^3\right\|_2^2
 +10Q-8V-16I_3-5I_4\\
 &\quad+\frac V{2\pi}(6V+4I_3+I_4-2Q).
 \end{aligned}
\end{equation}

\emph{Fourier closure.}  With $T=e^{is+i\varphi(s)+i\alpha}$, $\varphi_s=u$ and $\int\varphi=0$, the closure $\int T=0$ and $|e^{ix}-1-ix|\le x^2/2$ imply
\[
 \left|\int_0^{2\pi}e^{is}\varphi(s)ds\right|
 \le\frac12\norm\varphi_2^2\le\frac12V.
\]
Let $V_1$ be the squared $L^2$ norm of the $j=\pm1$ Fourier modes of $u$.  Parseval gives
\begin{equation}\label{eq:resonance}
 V_1\le\frac{V^2}{4\pi},\qquad
 \frac{X^2}{V}\ge16-\frac{15V}{4\pi},\qquad
 \frac QV\ge4-\frac{3V}{4\pi}.
\end{equation}
For $V>0$, define $z=X/\sqrt V$.  Since $d\le8/3$ and $\pi>157/50$, formula~\eqref{eq:resonance} gives $z>492/125$.
Modewise Fourier estimates yield
\begin{equation}\label{eq:Linfty}
 \norm u_\infty\le C_2\sqrt Q+W,
 \qquad \norm u_\infty\le C_4 X+W,
 \qquad W:=\frac V{2\pi},
\end{equation}
where
\[
 C_2=\sqrt{\frac\pi6-\frac1\pi}<\frac{227}{500},
 \qquad C_4=\sqrt{\frac{\pi^4/45-2}{2\pi}}<\frac{81}{500}.
\]
Indeed, the first mode satisfies $\norm{P_1u}_\infty\le W$ by~\eqref{eq:resonance}; for $|j|\ge2$ apply Cauchy--Schwarz to the Fourier series, using respectively
$\sum_{|j|\ge2}j^{-2}=\pi^2/3-2$ and
$\sum_{|j|\ge2}j^{-4}=\pi^4/45-2$.
The stated bounds on $C_2,C_4$ follow from $\pi<355/113$.

\emph{Retaining the quartic sign.}  Put $M_u=\norm u_\infty$, $s_4=\sqrt{I_4}/X$, and $a=(3+M_u)/2$.  Cauchy--Schwarz gives
\begin{equation}\label{eq:s-bounds}
 0\le s_4\le\frac{M_u}{z},\quad
 \frac{|I_3|}{X^2}\le\frac{s_4}{z},\quad
 \left\|\frac32u^2+\frac12u^3\right\|_2\le aXs_4.
\end{equation}
Fourier series of $u$ (which has zero mean) also give the elementary inequality
\begin{equation}\label{eq:quadratic-Fourier}
 10Q-8V\le2X^2,
\end{equation}
indeed $2j^4-10j^2+8=2(j^2-1)(j^2-4)\ge0$ for every nonzero integer $j$.
When $as_4<1$, using~\eqref{eq:s-bounds} and discarding the nonpositive term $-2WQ$ in~\eqref{eq:signed-evolution} gives
\begin{equation}\label{eq:f}
 \frac{U'}{LX^2}\le
 -4(1-as_4)^2+2+\frac{16+4W}{z}s_4+(-5+W)s_4^2+\frac{6W}{z^2}.
\end{equation}
We now bound its right side in three ranges.  In all of them,
\[
 x:=\sqrt V<x_*:=\frac{163}{250},\qquad
 W<w_0:=\frac{17}{250},\qquad z>z_0:=\frac{492}{125},
\]
as follows from $V\le4/(3\pi)$ and $\pi>157/50$.

If $z_0<z\le7$, the second estimate in~\eqref{eq:Linfty} and $C_4<81/500$ give $M_u\le c_*xz+W$, where $c_*:=81/500$.  Set
\[
 s_*(z)=c_*x_*+\frac{w_0}{z},\qquad
 a_*(z)=\frac12(3+c_*x_*z+w_0),\qquad
 r_*(z)=a_*(z)s_*(z).
\]
Then $s_4\le s_*(z)$ and $a\le a_*(z)$.  The function $r_*(z)$ is convex, so its maximum on either of the intervals below is attained at an endpoint.  Rational evaluation there gives
\[
\begin{array}{c|cc}
\text{range of }z & \text{bound on }r_*(z)&\text{bounds on }s_*(z)\\ \hline
 [z_0,5] & r_*<134/625 & 149/1250<s_*<12291/100000\\
 \lbrack5,7\rbrack & r_*<11/50 & 23/200<s_*<3/25
\end{array}
\]
In particular $r_*<11/50$ and $s_*<123/1000$.  The right side of~\eqref{eq:f} increases in $a,s_4,W$ here: its $s_4$ derivative is at least $12(1-11/50)-10(123/1000)>0$, and the other derivatives are nonnegative.  Substituting the bounds of the two rows into~\eqref{eq:f} gives, respectively,
\begin{equation}\label{eq:low-z}
 \frac{U'}{LX^2}<
 -\frac{67489765177}{15759375000000}< -\frac1{250},
 \qquad
 \frac{U'}{LX^2}<-\frac{919777}{10000000}<-\frac1{250}.
\end{equation}
For example, the first fraction is the exact value of
\[
 -4\left(1-\frac{134}{625}\right)^2+2
 +\frac{(16+4w_0)(12291/100000)}{z_0}
 +\frac{6w_0}{z_0^2}
 -(5-w_0)\left(\frac{149}{1250}\right)^2.
\]

If $z\ge7$, then $Q\le Vz$ and the first estimate in~\eqref{eq:Linfty} imply
\[
 M_u\le\left(C_2+\frac{\sqrt V}{2\pi\sqrt z}\right)\sqrt{Vz}
 <\left(\frac{227}{500}+\frac1{25}\right)\sqrt{Vz}
 <\frac12\sqrt{Vz}.
\]
Here $\sqrt V/(2\pi\sqrt z)<1/25$ follows from $x<x_*$, $z\ge7$, and $\pi>157/50$.
Thus $s_4\le\frac12\sqrt{V/z}$ and $a\le\widehat a:=(3+\frac12\sqrt{Vz})/2$.  Write $y:=z^{-1/2}$, so that $0<y\le7^{-1/2}<y_*:=189/500$.  Then
$\widehat a s_4\le\frac34xy+\frac18x^2<1/4$ and $s_4\le\frac12xy<1/8$, so, exactly as in the previous range, the right side of~\eqref{eq:f} increases with $a,s_4,W$.  We may therefore set $s_4=\frac12xy$, $a=\widehat a$, and $W=w_0$, and we obtain the polynomial bound
\begin{equation}\label{eq:poly}
 \begin{split}
 \frac{U'}{LX^2}\le\mathcal F(x,y):={}&-2+6xy+x^2+(8+2w_0)xy^3\\
 &-\frac{14-w_0}{4}x^2y^2-\frac34x^3y-\frac1{16}x^4+6w_0y^4.
 \end{split}
\end{equation}
The polynomial $\mathcal F$ is increasing in both $x$ and $y$ on the rectangle $(0,x_*)\times(0,y_*)$.  Indeed,
\begin{align*}
 \partial_x\mathcal F
 &=y\left(6-\frac{14-w_0}{2}xy-\frac94x^2\right)
 +x\left(2-\frac14x^2\right)+(8+2w_0)y^3>0,\\
 \partial_y\mathcal F
 &=x\left(6-\frac{14-w_0}{2}xy-\frac34x^2\right)
 +3(8+2w_0)xy^2+24w_0y^3>0
\end{align*}
there, since $xy<0.247$ and $x^2<0.426$.  Thus
\begin{equation}\label{eq:rational}
 \mathcal F(x,y)\le\mathcal F\left(\frac{163}{250},\frac{189}{500}\right)
 =-\frac{1605550457639}{15625000000000}< -\frac1{250}.
\end{equation}
This proves~\eqref{eq:diss-bound} when $V>0$.

\emph{Conclusion.} If $V=0$, the curvature is constant, the curve is a circle,
and $U'=0$. Thus, along the original free flow, $U'\le0$ whenever
$U-4\pi^2\le8/3$, with strict inequality for noncircular curves. The
dilation used above is frozen at the time of evaluation; it does not modify
the evolution.

Set $C=4\pi^2+8/3$ and assume $U(0)\le C$. Letting $t^*$ be the supremum of
the times $t$ for which $U(\tau)\le C$ for every $\tau\in[0,t]$, 
continuity and maximality imply $U(t^*)=C$. The curve at
$t^*$ is not a circle, so $U'(t^*)<0$. Hence $U<C$ on a right
neighbourhood of $t^*$, contradicting its definition. The sublevel is
therefore forward invariant, and $U$ is nonincreasing on it.
\end{proof}

\section{Centroid mobility}\label{sec:geom}
We retain the explicit dependence on the energy defect
\begin{equation}\label{eq:ad}
 d=U-4\pi^2,\qquad a(d)=\sqrt{\frac14+\frac{3d}{32\pi^2}},\qquad
 A(d)=\frac14+\frac{a(d)}\pi.
\end{equation}

\begin{lemma}\label{lem:angle}
Let $\gamma$ be a smooth closed unit-speed planar curve of turning number $+1$.  Write
\begin{equation}\label{eq:tangent-angle}
 k=\kappa+u,\quad\kappa=\frac{2\pi}L,\quad\int_0^Lu\,ds=0,\quad
 T=e^{i\alpha}e^{i\kappa s}e^{i\varphi(s)},\quad\varphi_s=u,\quad\int_0^L\varphi=0.
\end{equation}
With $V=\int_0^Lu^2\,ds$ and $d=LV$, the $j=\pm1$ curvature Fourier component satisfies
\begin{equation}\label{eq:V1}
 V_1\le\frac{L}{8\pi^2}V^2,
\end{equation}
and
\begin{equation}\label{eq:phase-bound}
 \norm\varphi_2\le\frac{L}{2\pi}\sqrt V\,a(d),\qquad
 \frac{\norm\varphi_2}{\sqrt L}\le\frac{\sqrt d}{2\pi}a(d).
\end{equation}
\end{lemma}
\begin{proof}
Closure gives $0=\int_0^Le^{i\kappa s+i\varphi(s)}ds$.  Using $\int_0^Le^{i\kappa s}ds=0$, $|e^{ix}-1-ix|\le x^2/2$ and Wirtinger's inequality $\norm\varphi_2\le\frac L{2\pi}\norm{\varphi_s}_2$,
\[
 \left|\int_0^Le^{i\kappa s}\varphi(s)ds\right|
 \le\frac12\norm\varphi_2^2\le\frac{L^2}{8\pi^2}V.
\]
Writing $\varphi=\sum_{j\ne0}\varphi_je^{i\kappa js}$ and $u=\sum_{j\ne0}u_je^{i\kappa js}$, so that $u_j=i\kappa j\varphi_j$ and $\int_0^Le^{i\kappa s}\varphi\,ds=L\varphi_{-1}$, this gives $|u_{-1}|\le V/(4\pi)$ and $V_1=2L|u_{-1}|^2\le LV^2/(8\pi^2)$.  Finally, by Parseval's identity and $j^2\ge4$ for $|j|\ge2$,
\[
 \norm\varphi_2^2\le\left(\frac L{2\pi}\right)^2
 \left(V_1+\frac14(V-V_1)\right)
 \le\left(\frac L{2\pi}\right)^2V\left(\frac14+\frac{3d}{32\pi^2}\right).
\]
This is~\eqref{eq:phase-bound}.
\end{proof}

\begin{lemma}\label{lem:mobility}
Let $H=k(\gamma-c)+\nu$.  Under the assumptions of Lemma~\ref{lem:angle},
\begin{equation}\label{eq:H-bound}
 \norm H_2\le A(d)L\sqrt V.
\end{equation}
Along the free elastic flow one consequently has
\begin{equation}\label{eq:b}
|c'|\le b(d)R_E',\qquad
 b(d)=\frac{A(d)(4\pi^2+d)\sqrt d}{4\pi^2}.
 \end{equation}
\end{lemma}
\begin{proof}
Set $T_0=e^{i\alpha}e^{i\kappa s}$, $\nu_0=iT_0$, and $\Gamma_0=-\nu_0/\kappa$.  Then $\Gamma_0'=T_0$, and $\int_0^L\Gamma_0\,ds=0$.  With $w=(\gamma-c)-\Gamma_0$, the function $w$ has zero mean and $w_s=T-T_0$.  By Wirtinger's inequality and $|e^{ix}-1|\le|x|$,
\[
 \norm w_2\le\frac L{2\pi}\norm{T-T_0}_2\le\frac L{2\pi}\norm\varphi_2,
 \qquad\norm{\nu-\nu_0}_2\le\norm\varphi_2.
\]
Since $\nu_0=-\kappa\Gamma_0$, we have $H=u(\gamma-c)+\kappa w+(\nu-\nu_0)$, and the centered bound~\eqref{eq:cent-bound} together with~\eqref{eq:phase-bound} gives
\[
 \norm H_2\le\frac L4\sqrt V+2\norm\varphi_2
 \le A(d)L\sqrt V.
\]
By~\eqref{eq:translation} and Cauchy--Schwarz, $|c'|\le L^{-1}\norm G_2\norm H_2\le A(d)\sqrt{d/L}\sqrt Y$.  Since $R_E'=2\pi Y/B^2$ and \eqref{eq:dilation} gives $\sqrt Y\ge2\pi B/L^{3/2}$,
\[
 \frac{|c'|}{R_E'}\le
 \frac{A(d)B^2\sqrt d}{2\pi\sqrt{LY}}
 \le\frac{A(d)BL\sqrt d}{4\pi^2}
 =\frac{A(d)(4\pi^2+d)\sqrt d}{4\pi^2}.
\]
\end{proof}

\begin{lemma}\label{lem:margin}
Under the assumptions of Lemma~\ref{lem:angle}, let
$S_*:=\pi^2/3-1$.  Then for every unit vector $e$,
\begin{equation}\label{eq:margin-bound}
 \rho_e^\pm\ge\frac L{2\pi}\left[
 1-\frac{d\,a(d)^2}{8\pi^2}
 -\sqrt{S_*}\,\frac{\sqrt d\,a(d)}{2\pi}\right].
\end{equation}
\end{lemma}
\begin{proof}
Identify $\R^2$ with $\mathbb C$ and expand $\gamma-c=\sum_{j\ne0}\hat\gamma_je^{i\kappa js}$ and $T=\sum_j\hat T_je^{i\kappa js}$ in complex Fourier series, so that $\hat\gamma_j=\hat T_j/(i\kappa j)$ for $j\ne0$.  The coefficient $\hat\gamma_1$ satisfies
\[
 |\hat\gamma_1|=\frac L{2\pi}\left|\frac1L\int_0^Le^{i\varphi}\,ds\right|
 \ge\frac L{2\pi}\left(1-\frac{\norm\varphi_2^2}{2L}\right),
\]
using $\cos x\ge1-x^2/2$.  For $j\ne0,1$, the coefficient $\hat T_j$ coincides with the $j$-th coefficient of $T-T_0$, since $T_0$ is a pure first harmonic.  By Cauchy--Schwarz and Parseval,
\[
 \norm{\sum_{j\ne0,1}\hat\gamma_je^{i\kappa js}}_\infty
 \le\frac L{2\pi}\sqrt{\sum_{j\ne0,1}\frac1{j^2}}
 \frac{\norm{T-T_0}_2}{\sqrt L}
 \le\frac L{2\pi}\sqrt{S_*}\frac{\norm\varphi_2}{\sqrt L}.
\]
The leading harmonic $\hat\gamma_1e^{i\kappa s}$ describes a circle of radius $|\hat\gamma_1|$ centered at the origin.  In each direction its extreme values are $\pm|\hat\gamma_1|$; the remainder changes either extremum by at most its $L^\infty$ norm.  Substituting~\eqref{eq:phase-bound} proves~\eqref{eq:margin-bound}.
\end{proof}

\begin{lemma}\label{lem:numerical}
For every $0\le d\le8/3$,
\begin{equation}\label{eq:constants}
 b(d)<\frac34,\qquad
 1-\frac{d\,a(d)^2}{8\pi^2}
 -\sqrt{\frac{\pi^2}{3}-1}\frac{\sqrt d\,a(d)}{2\pi}
 >\frac34.
\end{equation}
\end{lemma}
\begin{proof}
The first expression increases and the second decreases with $d$, so it suffices to verify the endpoint $d=8/3$.  The elementary bounds
\[
 \frac{157}{50}<\pi<\frac{22}7,\quad
 \sqrt{\frac83}<\frac53,\quad
 a(8/3)<\frac{53}{100},\quad
 \sqrt{\frac{\pi^2}{3}-1}<\frac{303}{200},\quad
 4\pi^2+\frac83<\frac{169}{4}
\]
lead to
\[
 b(d)<\frac{(21/50)(169/4)(5/3)}{4(157/50)^2}
 =\frac{147875}{197192}<\frac34.
\]
Similarly,
\[
 \begin{split}
 &1-\frac{d\,a(d)^2}{8\pi^2}
 -\sqrt{\frac{\pi^2}{3}-1}\frac{\sqrt d\,a(d)}{2\pi}\\
 &\qquad>1-\frac{(8/3)(53/100)^2}{8(157/50)^2}
 -\frac{(303/200)(5/3)(53/100)}{2(157/50)}
 =\frac{9197897}{11831520}>\frac34.
 \end{split}
\]
\end{proof}

\section{The length-preserving flow}\label{sec:lp}
In this section $G=G_c$ and $\lambda_c$ is given by~\eqref{eq:lambda-intro}. Thus $G_c$ is the $L^2(ds)$-orthogonal projection of $G_0$ onto the orthogonal complement of $k$, which is the tangent space of the fixed-length constraint; the multiplier $\lambda_c$ is not assumed to be positive. Set
\begin{equation}\label{eq:constrained-Y}
 Y_c=\int_\gamma G_c^2\,ds.
\end{equation}

\begin{lemma}\label{lem:constrained-variations}
Along the length-preserving flow,
\begin{equation}\label{eq:constrained-dissipation}
 L'=0,\qquad B'=-Y_c,\qquad U'=-LY_c,
\end{equation}
and
\begin{equation}\label{eq:constrained-centroid}
 Lc'=\int_\gamma G_c\,[k(\gamma-c)+\nu]ds.
\end{equation}
\end{lemma}
\begin{proof}
The first variation of $B$ is $\int G_0f\,ds$, whereas that of $L$ is $-\int kf\,ds$. By definition $\int kG_c\,ds=0$. Thus, for $f=-G_c$, the length is fixed and
\[
 B'=-\int G_0G_c\,ds
 =-\int G_c^2\,ds-\lambda_c\int kG_c\,ds=-Y_c.
\]
Translation invariance of $B$ and $L$ gives $\int G_0\nu\,ds=\int k\nu\,ds=0$; hence also $\int G_c\nu\,ds=0$. Differentiating the centroid proves~\eqref{eq:constrained-centroid}, exactly as in Lemma~\ref{lem:exact}.
\end{proof}

\begin{proposition}\label{prop:constrained-gap}
Suppose that $\gamma$ has turning number $\pm1$ and $0<d:=LB-4\pi^2\le8/3$. Then
\begin{equation}\label{eq:constrained-gap}
 Y_c\ge\frac{25}{4}\frac{(2\pi)^4}{L^5}\,d.
\end{equation}
In particular, a noncircular curve in this sublevel cannot be a stationary point of the length-preserving flow.
\end{proposition}
\begin{proof}
By reversal of orientation if needed, we assume turning number $+1$. Normalize the \emph{curve at the time under consideration} to $L=2\pi$ and write $k=1+u$, $\int u\,ds=0$. Put
\[
 V=\int u^2\,ds=\frac d{2\pi},\qquad Q=\int u_s^2\,ds,
 \qquad I_j=\int u^j\,ds\quad(j=3,4).
\]
Since $d\le8/3$, we have $0<V\le4/(3\pi)<1/2$. Formula~\eqref{eq:resonance}, which depends only on closure and not on the evolution, yields
\begin{equation}\label{eq:constrained-q}
 q:=\frac QV\ge4-\frac{3V}{4\pi}>\frac{31}{8}.
\end{equation}
Since $B=2\pi+V$ and
\[
 \lambda_c=\frac{2\pi+6V+4I_3+I_4-2Q}{2\pi+V},\qquad
 -\int G_cu\,ds=2Q-3V-3I_3-I_4+\lambda_c V,
\]
we obtain the exact identity
\begin{align}
 -\int G_cu\,ds={}&\frac{4\pi}{2\pi+V}(Q-V)
 +\frac{3V^2}{2\pi+V}\nonumber\\
 &+\left(\frac{4V}{2\pi+V}-3\right)I_3
 -\frac{2\pi}{2\pi+V}I_4.
 \label{eq:constrained-identity}
\end{align}
The first Fourier bound in~\eqref{eq:Linfty} is geometric and applies to the constrained flow. Together with $Q=Vq$, it gives
\begin{equation}\label{eq:constrained-uinfty}
 M:=\norm u_\infty\le C_2\sqrt{Vq}+\frac V{2\pi}
 <m(q):=\frac{\sqrt q}{3}+\frac{2}{25}.
\end{equation}
Indeed, since $V\le4/(3\pi)$ and $3<\pi<22/7$, we have $C_2^2V\le 2/9-4/(3\pi^2)<1/9$ and $V/(2\pi)\le2/(3\pi^2)<2/25$.  Now $|I_3|\le MV$ and $I_4\le M^2V$. Since $d\le8/3$ and $\pi>3$, we have $4\pi/(2\pi+V)>24/13$; also the coefficient of $I_3$ in~\eqref{eq:constrained-identity} has absolute value at most $3$. Discarding its positive term $3V^2/(2\pi+V)$, we obtain
\begin{equation}\label{eq:f-q}
 -\frac1V\int G_cu\,ds\ge h(q)
 :=\frac{24}{13}(q-1)-3m(q)-m(q)^2.
\end{equation}
The function $h$ is strictly increasing on $[31/8,\infty)$ since
\[
 h'(q)=\frac{24}{13}-\frac19-\frac{79}{150\sqrt q}>0.
\]
Because $\sqrt{31/8}<2$, we have $m(31/8)<56/75$. Hence
\begin{equation}\label{eq:f-lower}
 h(q)>h(31/8)>
 \frac{24}{13}\frac{23}{8}-3\frac{56}{75}
 -\left(\frac{56}{75}\right)^2
 =\frac{183557}{73125}>\frac52.
\end{equation}
By Cauchy--Schwarz,
\[
 \frac52V\le-\int G_cu\,ds
 \le \sqrt{Y_cV},
 \qquad\text{so}\qquad Y_c\ge\frac{25}{4}V.
\]
Finally, under a dilation of factor $2\pi/L$, $G_c$ scales as length$^{-3}$ and $Y_c$ as length$^{-5}$. Since $V=d/(2\pi)$ in the normalized coordinates, this becomes~\eqref{eq:constrained-gap} at arbitrary length.
\end{proof}

\begin{proposition}\label{prop:constrained-motion}
Under the hypothesis of Theorem~\ref{th:main} for the length-preserving flow, set $A_*=A(8/3)$, with $A$ defined in~\eqref{eq:ad}. Then $A_*<21/50$ and, for all $t_0<t_1$,
\begin{equation}\label{eq:constrained-displacement}
 \frac{|c(t_1)-c(t_0)|}{L}
 \le\frac{A_*}{10\pi^2}\bigl(U(t_0)-U(t_1)\bigr)
 <\frac1{8\pi^2}.
\end{equation}
Moreover, for every unit vector $e$ and every $t\in[0,+\infty)$,
\begin{equation}\label{eq:constrained-margins}
 \rho_e^\pm\ge \frac34\frac L{2\pi}.
\end{equation}
\end{proposition}
\begin{proof}
The energy is nonincreasing by~\eqref{eq:constrained-dissipation}; hence $d(t)\le8/3$, and Lemma~\ref{lem:topology} ensures turning number $\pm1$. The $L^2$ estimate for $H=k(\gamma-c)+\nu$ in~\eqref{eq:H-bound} is geometric and applies unchanged. By~\eqref{eq:constrained-centroid} and Cauchy--Schwarz,
\[
 |c'|\le A(d)\sqrt{\frac dL}\sqrt{Y_c}.
\]
For $d>0$, Proposition~\ref{prop:constrained-gap} and $-U'=LY_c$ imply
\begin{equation}\label{eq:constrained-speed}
 \frac{|c'|}{L}\le\frac{A(d)}{10\pi^2}(-U')
 \le\frac{A_*}{10\pi^2}(-U').
\end{equation}
At $d=0$ the curve is a stationary circle, so the same inequality holds. The elementary bounds in the proof of Lemma~\ref{lem:numerical} give $a(8/3)<53/100$ and therefore $A_*<1/4+53/(100\pi)<21/50$. Integrating~\eqref{eq:constrained-speed} on $[t_0,t_1]$ and using $0\le U(t_0)-U(t_1)\le8/3$ proves~\eqref{eq:constrained-displacement}, since $\frac{21}{50}\cdot\frac{8}{3}\cdot\frac1{10}=\frac{14}{125}<\frac18$. The geometric margin estimate~\eqref{eq:margin-bound}, together with Lemma~\ref{lem:numerical}, gives~\eqref{eq:constrained-margins}.
\end{proof}

\section{Proof of the main result}\label{sec:proof}
\begin{proof}[Proof of Theorem~\ref{th:main}]
Lemma~\ref{lem:topology} shows that the initial turning number is $\pm1$, which is preserved throughout the smooth evolution. We reverse orientation if necessary and assume it is $+1$.

For the free flow, Proposition~\ref{prop:persist} gives~\eqref{eq:main-persist}; the lower bound $U\ge4\pi^2$ is~\eqref{eq:oscillation} with $n=\pm1$. The estimates of Lemmas~\ref{lem:mobility}, \ref{lem:margin} and~\ref{lem:numerical}, together with $L/(2\pi)\ge R_E$, give~\eqref{eq:barrier-intro}. If~\eqref{eq:halfplanes} held, then
\[
 (c(t_0)-c(t_1))\cdot e
 \ge\rho_e^-(t_0)+\rho_e^+(t_1)
 \ge\frac34[R_E(t_0)+R_E(t_1)].
\]
On the other hand, integrating the bound $|c'|\le\frac34R_E'$ of~\eqref{eq:barrier-intro} and using $R_E'>0$, we get
\[
 |c(t_1)-c(t_0)|\le\frac34
 [R_E(t_1)-R_E(t_0)],
\]
Since $|(c(t_0)-c(t_1))\cdot e|\le|c(t_1)-c(t_0)|$, the two inequalities give $R_E(t_0)\le0$, which is impossible.

For the length-preserving flow, monotonicity~\eqref{eq:main-persist} follows directly from~\eqref{eq:constrained-dissipation}. If the two half-plane conditions held,~\eqref{eq:constrained-margins} would give
\[
 (c(t_0)-c(t_1))\cdot e
 \ge\rho_e^-(t_0)+\rho_e^+(t_1)
 \ge\frac{3L}{4\pi}.
\]
This contradicts~\eqref{eq:constrained-displacement}, as $L/(8\pi^2)<3L/(4\pi)$. In both cases the persistent bound $d(t)\le8/3$ and the second assertion of Lemma~\ref{lem:topology} also give embeddedness at every time. This proves the theorem.
\end{proof}

\begin{corollary}\label{cor:barrier}
Under the free-flow hypotheses of Theorem~\ref{th:main}, for every $e\in S^1$ both
\[
 \frac34R_E(t)+c(t)\cdot e,\qquad
 \frac34R_E(t)-c(t)\cdot e
\]
are nondecreasing. Moreover $\rho_e^\pm\ge(3/4)R_E$.
\end{corollary}
\begin{proof}
The derivatives of these two functions are $\frac34R_E'\pm c'\cdot e\ge0$ by~\eqref{eq:barrier-intro}.
\end{proof}

\begin{remark}
Our result concerns the free and length-preserving equations in the stated
energy sublevel. Huisken's original question is about the distinct
length-penalized equation~\cite{MPPSurvey,KemmochiMiura}.
\end{remark}

\section*{Acknowledgments}
The author is a member of the INdAM-GNAMPA research group.

\medskip
\noindent\textbf{Use of AI.}
The proofs in this paper were found with the assistance of a large language
model (LLM), which was also used to draft the text. The author has verified
all results and the text, and takes full responsibility for the content.

\end{document}